\documentclass[3p]{elsarticle}

\usepackage{graphicx} 
\usepackage{subfigure} 
\usepackage{upgreek}
\usepackage{capt-of}
\usepackage[usenames,dvipsnames]{color}
\usepackage{amsmath}
\usepackage{amssymb}

\newcommand{\dif}{\textup d}
\newcommand{\Real}{\mathbb R}
\newcommand{\Nat}{\mathbb N}

\usepackage[normalem]{ulem}

\journal{...}

\begin{document}

\title{Recent results in the systematic derivation and convergence of SPH}
\tnotetext[t1]{This work was awarded the \emph{Libersky prize} for the best student paper and presentation at the 11th SPHERIC International Workshop, 13-16 June 2016, Munich -- Germany.}

\author[1]{Iason Zisis}
\author[2]{Joep H.M. Evers}
\author[1]{Bas van der Linden}
\author[3]{Manh Hong Duong}

\address[1]{Department of Mathematics \& Computer Science, TU Eindhoven, The Netherlands}
\address[2]{Department of Mathematics \& Statistics, Dalhousie University, Halifax, Canada}
\address[3]{Mathematics Institute, University of Warwick, United Kingdom}

\cortext[cor1]{I. Zisis is the corresponding author (iason.zisis@outlook.com). J.H.M. Evers is supported by an AARMS Postdoctoral Fellowship and M.H. Duong is supported by ERC Starting Grant 335120.}
 
\maketitle

\section{Introduction}

In the literature of the Smoothed Particle Hydrodynamics (SPH) computational method, it is established that the standard SPH scheme can be derived by applying the principle of least action to a particle system, where the SPH density estimate acts as a constraint \cite{Bonet1999,Monaghan2005,Price2012}. Nonetheless, a subtlety lies in the fact that in the derivation of the SPH equations \emph{the action of the particle system is minimized rather than the action of the continuum}. 
The procedure of starting from the continuum, minimizing its action and then discretizing the equations reveals the mathematical formalities necessary to convince oneself that SPH indeed comes from principles of continuum mechanics.

The theory of \emph{measure-valued evolution equations} was first shown by Di Lisio \emph{et al.} \cite{DiLisio1998} to provide the adequate framework for the study of both the SPH particle system and the limiting continuum setting in a single context. Following them, in Evers \emph{et al.}  \cite{Evers2015} we adopt \emph{the Wasserstein distance in the space of probability measures}, which determines how ``close" two measures are by computing the (optimal) cost of transforming one into the other. 
Its definition allows to conveniently derive suitable upper bounds on the distance between measures, which are necessary ingredients to prove convergence. 
After constructing a discrete approximation of the initial measure, we prove the convergence of measure-valued solutions \cite{Evers2015}. There are two major differences with Di Lisio \emph{et al.} \cite{DiLisio1998} (and accordingly with a recent review \cite{Colagrossi2014}). First, their scheme is not the classic SPH scheme, but rather one that is known not to conserve momentum, whenever applied to relevant physical processes. Our proof applies to both the traditional SPH scheme and the former one. Second, we allow for a much more general class of force fields, including external and internal conservative forces, as well as friction and non-local interactions.
%
%

The overall aim of the present paper is to report recent developments in the formal derivation of SPH ---focusing on the authors' work \cite{Evers2015}--- and eventually bridge the gap between the SPH literature and the few publications dealing with the rigorous mathematical formalism of SPH. The structure of the paper is the following:
Section \ref{sec:SPH from cm} summarizes the systematic procedure followed for the derivation of measure-valued and particle formulations of continuum mechanics equations. 
Section \ref{sec: Conv} states the theoretical convergence result. Section \ref{sec: add} discusses on the limitations of the proof. Section \ref{sec: numerics} presents numerical paradigms which exhibit the applicability of the convergence with respect to the Wasserstein distance.
%
Finally, the related open problems are stated in Section \ref{sec: concl}.

\section{SPH from continuum mechanics}\label{sec:SPH from cm}

Denote the particle trajectories $x = x(x_0,t) \in \Omega_t \subset \Real^d$, taken with respect to an initial configuration of the medium $x_0 \in \Omega_0 \subset \Real^d$, with $d$ being the spatial dimension. In the following, on grounds of concise notation, we suppress the dependence on $t$ and it is only implied. For a medium found in the domain $\Omega$, its mass:
\begin{equation}\label{eq: mass_measure}
  \mu(\Omega) = \int_{\Omega} \dif \mu(x) = \int_{\Omega_0} \dif \mu_0(x_0) = \mu_0(\Omega_0),
\end{equation}
is always conserved in the absence of sinks or sources. The above equality is guaranteed by the one-to-one correspondence between the set of particles in the medium's reference configuration $\Omega_0$ and any later configuration $\Omega$. Thus, for the integral of some bounded and measurable $f$ with respect to the measure $\mu$, it is possible to perform a coordinate transform from $\Omega$ to $\Omega_0$ and obtain 
$
\int_\Omega f(x) \, \dif \mu(x) =	\int_{\Omega_0} f(x) \, \dif \mu_0(x_0),
$
where we recall that $x=x(x_0,t)$. 
%
In the context of measures the \emph{mass density} function is defined as the (\emph{Radon-Nikodym}) derivative of $\mu$ with respect to the \emph{Lebesgue measure} $\lambda$:
\begin{equation}\label{eq: density}
  \rho = \frac{\dif \mu}{\dif \lambda}.
\end{equation}
Consequently, (\ref{eq: mass_measure}) attains the typical form:
$
  \int_{\Omega} \rho(x) \, \dif\lambda(x) = \int_{\Omega_0} \rho_0(x_0) \, \dif\lambda(x_0).
$
The measure-valued definition (\ref{eq: mass_measure}) of mass conservation is more general, since it does not require $\rho$ to be well-defined. Thus, it offers the basis for a generalized approach for both continuous and discretized media. Within this context, we may construct a discrete approximation $\mu_0^N$ of the measure of mass $\mu_0$, which then can be used for numerical analysis. The one-to-one correspondence between particles in $\Omega_0$ and $\Omega$ and the definition of the medium's density function (\ref{eq: density}) deliver 
$
	\rho(x) \, \dif\lambda(x) = \rho_0(x_0) \, \dif\lambda(x_0),
$
%
%
and consequently the medium's density is given by:
\begin{equation}\label{eq: density total}
	\rho(x) = \frac{\rho_0(x_0)}{J(x)},
\end{equation}
where $J$ is the determinant of the Jacobian matrix \cite[Eq. (5)]{Seliger1968}.

The action of the medium,
$
\mathcal S = \int_{0}^{T} L \, \dif t,
$
involves the Lagrangian: 
\begin{equation*}
L = \int_\Omega \Big ( \frac{1}{2} \|\dot x\|^2 - e(\rho(x)) \Big ) \rho(x)  \, \dif\lambda(x),
\end{equation*}
where $e$ is the internal energy of the medium prescribed by the thermodynamic relation $\dif e/\dif \rho = P(\rho)/\rho^2$ for ideal processes, with $P$ the mean volumetric stress --- \emph{pressure}. According to the principle of least action, the equations of motion follow from the minimization problem
$
 \updelta \mathcal S |_{0}^{T} = \int_{0}^{T} \updelta L \, \dif t = 0,
$
where the differential $\updelta$ denotes a variation of the particle trajectories $x$. Therefore, in order to obtain the variation and eventually the equations of motion, the spatial domain of integration (in the Lagrangian) should be independent of the variation. This is achieved either by writing the Lagrangian with respect to the initial configuration or in the context of measures.

In particular, three steps are necessary to pass from the action of the continuous system to the motion equations of the resulting particle system:
\begin{itemize}
\item{A:} Introduce the measure-valued formulation by replacing $\rho \, \dif\lambda$ with $\dif\mu$ and, wherever necessary, approximate $\rho$ by some $\tilde \rho$, which depends on the measure $\mu$. Typically, this is:
\begin{equation}\label{eq: mass_conservation}
	\tilde \rho(x) = \int_\Omega W_h(x-y) \dif\mu(y),
\end{equation}
where $W_h$ is a symmetric positive mollifier, see e.g., \cite{Monaghan2005}.
\item{B:} Substitute for $\mu$ the discrete measure:
\begin{equation}\label{eq: discrete}
	\mu^N = \sum_{i=1}^N m_i \, \delta_{x_i}.
\end{equation}
where $\delta_{x_i} = 1$ if $x = x_i$ and zero otherwise, is the Dirac measure at $x_i$.
\item{C:} Derive the Euler-Lagrange equations.
\end{itemize}
These three steps have been introduced in more details in \cite{Evers2015}. Step A takes us to a regularized version of the problem, which is a problem different from the original one. Step B cannot happen before A, but we have the freedom to choose the further ordering. This gives rise to three different derivations:
\begin{itemize}
\item{ABC:} Write the measure-valued Lagrangian:
\begin{equation*}
	\tilde L = \int_\Omega \Big ( \frac{1}{2} \|\dot x\|^2 - e(\tilde \rho(x)) \Big ) \dif\mu(x),
\end{equation*}
discretize it:
\begin{equation*}
	\tilde L^N= \sum_i \Big ( \frac{1}{2} \|\dot x_i^2\| - e(\tilde \rho_i) \Big ) m_i,
\end{equation*}
and derive the corresponding equations of motion afterwards:
\begin{equation*}
	\frac{d}{dt}\Big ( \frac{\partial \tilde L^N}{\partial \dot x_i} \Big ) - \frac{\partial \tilde L^N}{\partial x_i} = 0,
\end{equation*}
\begin{equation}\label{eq: sph_motion1}
	\ddot x_i = -\sum_j \Big ( \frac{P(\tilde \rho_i)}{\tilde \rho_i^2} + \frac{P(\tilde \rho_j)}{\tilde \rho_j^2} \Big ) \, \nabla W_h(x_i-x_j) \, m_j.
\end{equation}
It is the common technique encountered in the SPH literature \cite{Bonet1999,Monaghan2005,Price2012} and its importance was recognized already in early articles about SPH (e.g. \cite{Monaghan1978}).

\item{ACB:} From the measure-valued Lagrangian: 
\begin{equation*}
	\tilde L = \int_\Omega \Big ( \frac{1}{2} \|\dot x\|^2 - e(\tilde \rho(x)) \Big ) \dif\mu(x),
\end{equation*}
derive the equations of motion:
\begin{equation*}
	\frac{d}{dt}\Big ( \frac{\partial \tilde L}{\partial \dot x} \Big ) - \frac{\partial \tilde L}{\partial x} = 0,
\end{equation*}
\begin{equation}\label{eq: continuous}
	\ddot x = -\int_\Omega \Big ( \frac{P(\tilde \rho)}{\tilde \rho^2} \Big |_{(x)} + \frac{P(\tilde \rho)}{\tilde \rho^2} \Big |_{(y)} \Big ) \, \nabla W_h(x-y) \, \dif\mu(y),
\end{equation}
and discretize these equations afterwards:
\begin{equation}\label{eq: sph_motion2}
	\ddot x_i = -\sum_j \Big ( \frac{P(\tilde \rho_i)}{\tilde \rho_i^2} + \frac{P(\tilde \rho_j)}{\tilde \rho_j^2} \Big ) \, \nabla W_h(x_i-x_j) \, m_j.
\end{equation}
It is the technique introduced in the authors' work \cite{Evers2015}.

\item{CAB:} From the Lagrangian of the continuous system in the reference configuration:
\begin{equation*}
	L = \int_{\Omega_0} \Big ( \frac{1}{2} \|\dot x\|^2 - e(\rho_0(x_0)/J(x)) \Big ) \rho_0(x_0)  \, \dif\lambda(x_0),
\end{equation*}
derive the equations of motion taking into account the functional relation $\rho(x) = \rho_0(x_0)/J(x)$ as in \cite{Seliger1968}:
\begin{equation*}
	\frac{d}{dt}\Big ( \frac{\partial L}{\partial \dot x} \Big ) - \frac{\partial L}{\partial x} = 0,
\end{equation*}
\begin{equation*}
	\ddot x = -\frac{1}{\rho} \frac{\dif P(\rho)}{\dif\rho} \nabla \rho,
\end{equation*}
write them in the measure-valued form:
\begin{equation*}
	\ddot x = -\frac{1}{\tilde \rho(x)} \frac{\dif P(\tilde \rho)}{\dif\tilde \rho} \Big |_{(x)} \int_\Omega \nabla W_h(x-y) \dif\mu(y),
\end{equation*}
and finally discretize them:
\begin{equation*}
	\ddot x_i = -\frac{1}{\tilde \rho_i} \frac{\dif P(\tilde \rho)}{\dif \tilde \rho} \Big |_i \sum_j \nabla W_h(x_i-x_j) m_j.
\end{equation*}
This strategy is implied by \cite{DiLisio1998,Colagrossi2014}.
\end{itemize}

It thus turns out that the order in which these steps are executed determines what the resulting equation is. To be more precise, the classical SPH scheme \cite{Price2012} is obtained, whenever the regularization of the density takes place before applying the principle of least action (A-B-C and A-C-B). If we apply the principle of least action to the action at the continuum level before regularizing the density (C-A-B), then we obtain the scheme appearing in \cite{DiLisio1998,Colagrossi2014}, which is rarely ever employed for SPH computations. 
We emphasize that although both schemes arrive from the principle of least action, the latter can also be derived directly from Newtonian mechanics and introduction of the density regularization. 
Finally, a parameter $\theta$ may be used to put the two equations in the common formulation:
\begin{eqnarray}\label{eq: system}
	\ddot x = - F_\theta(\tilde \rho(x)) \int_\Omega  \nabla W_h(x-y)\, \dif\mu(y) - \theta \int_\Omega F_\theta(\tilde \rho(y)) \nabla W_h(x-y)\, \dif\mu(y),
\end{eqnarray}
for $\theta \in \{0,1\}$, where $F_0(\tilde \rho) = 1/\tilde \rho \, d(P(\tilde \rho))/\dif\tilde \rho$ and $F_1(\tilde \rho) = P(\tilde \rho)/\tilde \rho^2$.

%
\section{Convergence result}
\label{sec: Conv}

The \emph{Wasserstein distance} between two probability measures $\mu_1$ and $\mu_2$ is defined as:
\begin{equation*}
	\mathcal W(\mu_1, \mu_2) = \inf_{\pi \in \Pi(\mu_1,\mu_2)} \int_{\Real^d \times \Real^d} |\chi-\psi| \, \pi(\dif\chi,\dif\psi),
\end{equation*}
where $\Pi(\mu_1,\mu_2)$ is the set of all \emph{joint representations} of $\mu_1$ and $\mu_2$. Joint representations are also called \textit{couplings} and are defined such that for each $i=1,2$, 
$
\int_{\Real^d\times \Real^d} f(\chi_i) \, \pi(\dif\chi_1,\dif\chi_2) = \int_{\Real^d} f(\chi) \, \dif\mu_i(\chi),
$
for all measurable, bounded functions $f$ on $\Real^d$. As said in the introduction, effectively, the Wasserstein distance computes the (optimal) cost of transforming one probability measure into another. An exposition on the Wasserstein distance and the related concept of \emph{optimal transport}, can be found in \cite{Villani2009}. Considering the system of (\ref{eq: system}) and (\ref{eq: mass_conservation}), written with respect to the absolutely continuous measure $\mu_t$, and the corresponding system written for the discrete measure $\mu_t^N$ from (\ref{eq: discrete}), the proof in Evers \emph{et al.} \cite{Evers2015} establishes that:
\begin{equation*}
	\sup_{t\in [0,T]}\mathcal{W}(\mu^N_t,\mu_t) \rightarrow 0, \,\,\, \textup{as} \,\,\,  N\rightarrow \infty,
\end{equation*}
provided that we can approximate the initial measure arbitrarily well. This holds for both, the classical SPH scheme and the non-conservative one of Di Lisio \emph{et al.} \cite{DiLisio1998} as well; $\theta = 1$ and $\theta=0$ respectively in (\ref{eq: system}). For clarity, here we explicitly write the  dependence of the measures on time. 

The paradigms of the numerical illustration of Section \ref{sec: numerics} are performed for a series of increasing particle numbers $N_k$. For each computation, after normalizing the total mass of the system, at each time instance $t\in I:= \{\varphi T/10$, $\varphi=\{0,...,10\}\subset \Nat\}$, we solve a linear programming problem to calculate: 
%
%
%
%
%
%
%
\begin{equation}\label{eqn: approx sup wass}
M_{k,k+1} := \max_{t\in I} \mathcal{W}(\mu^{N_k}_t,\mu^{N_{k+1}}_t) \approx \sup_{t\in [0,T]}\mathcal{W}(\mu^{N_k}_t,\mu^{N_{k+1}}_t), \hspace{0.1in} C_{k+1}^{(d)} := \log_{\frac{N_{k+1}}{N_k}} \Big ( \frac{M_{k+1,k+2}}{M_{k,k+1}} \Big ).
\end{equation}
%
%
%
The theoretically predicted convergence rate is the same as for the initial measure, i.e. $O(N^{-1/d})$, whence we expect that $C^{(d)}_{k+1}$ tends to the value $-1/d$.

A critical point of the theoretical result of Evers et al. \cite{Evers2015} (and Di Lisio \emph{et al.} \cite{DiLisio1998}) is that it makes no conclusion on the smoothing length $h$. The convergence proof is achieved for $h$ fixed with the number of particles, and the dependence of $h$ on $N$ is not investigated. It is known that in order for the regularized equations of hydrodynamics to approximate the real physics well, $h$ should be sufficiently small. In the SPH literature (e.g. \cite{Monaghan2005}), it is common practice to achieve this by taking $h = \eta \, N^{-1/d}$, with parameter $1.2\leq\eta\leq1.5$, for Gaussian-like kernels. By extension, cases of spatially and temporally varying $h$, like those used in shock problems \cite{Monaghan2005}, are not covered by the theoretical result.

\section{Additional processes and limitations of the theoretical result}
\label{sec: add}
The theoretical proof of convergence \cite{Evers2015} covers cases broader than the one discussed in Section \ref{sec:SPH from cm}. First and foremost, the form of the potential energy covered by the theoretical proof of convergence is
$
	e^* = e^*(\rho(x),x) = e(\rho(x)) + u(x),
$
with $e$ the internal energy of the medium given by $\partial e/\partial \rho = P(\rho)/\rho^2$, as discussed in Section \ref{sec:SPH from cm}, and $u$ an external field, such as gravity. A limitation of the theoretical proof is that the admissible equations of state are of the form
$
	P(\rho) = \mathcal{K}\rho^\gamma,
$
where $\mathcal{K}$ is a parameter and $\gamma$ is the so-called \textit{polytropic exponent}, which needs to satisfy $\gamma > 1$. It should be underlined that the theoretical proof is not conclusive about equations of state in the form $P(\rho) = B \cdot((\rho/\rho_0)^\gamma - 1)$, which are typically employed in SPH computations for the modeling of water \cite{Monaghan2005}.

Processes involving $e^*$ are conservative, and therefore, the related equation of motion follows the procedures of Section \ref{sec:SPH from cm}. Apart from them, the proof of convergence covers processes described by the equation:
$
	\updelta \mathcal S |_{0}^{T} = -\updelta \mathcal Q |_{0}^{T},
$
where:
\begin{equation*}
	\updelta \mathcal Q |_{0}^{T} = \int_{0}^{T} \int_\Omega \Big ( - \nu(x) \, \dot x + \int_\Omega K(x-y) \, \dif\mu(y) \Big ) \, \updelta x \,\dif\mu(x) \, \dif t.
\end{equation*}
In the general case $K$ can be an anisotropic kernel, $K(x-y) \neq K(\|x-y\|)$, describing non-local interactions within the system. On the other hand, $\nu = \nu(x)$ is a dissipative term, without non-local characteristics. It should be stressed that this is in contrast to the non-locality of the dissipative term:
\begin{equation}\label{eq: nonlocal_dissipation}
	\sim \int_\Omega (\dot x - \dot y) \nabla^2 W_h(x-y) \, \frac{\dif\mu}{\tilde \rho} \Big |_{(y)},
\end{equation}	
typically constructed (via approximation of $\nabla^2 W_h$) in SPH to model viscosity (e.g. \cite{Monaghan2005}). Note that this construction further assumes the approximation: $\dif\lambda = \dif\mu/\rho \approx \dif\mu/\tilde \rho$, which is necessary if in the SPH-approximation of the continuity equation \cite{Monaghan2005,Monaghan2013} the latter approximation is preferred to describe mass conservation (see Section \ref{sec: numerics}) over the temporal evolution of (\ref{eq: mass_conservation}).

\section{Numerical paradigms}\label{sec: numerics}

We construct the initial measure $\mu^N_0$, corresponding to the $N$-particle approximation of $\mu_0$, according to a partitioning of the initial domain in $N$ subdomains of incremental volume $V_i$. Masses are assigned as $m_i = \rho_0(x_i)\,V_i$ for each $i=1,\ldots,N$. In Evers \emph{et al.} \cite[Section 3.5]{Evers2015} we show two formal ways of constructing the sequence $\mu^N_0$ such that it converges to $\mu_0$ at rate $\mathcal{O}(N^{-1/d})$; they correspond to particle initialization strategies typically used in the SPH literature. The theoretical convergence result \cite{Evers2015} establishes that the corresponding solutions $\mu^N$ converge at the same rate.
%
%
%
In Evers \emph{et al.} \cite{Evers2015}, test cases which conform to the assumptions of the proof are examined. The cases that follow here, suggest that the same theoretical results may be expected to hold also for cases that are typically used for benchmarking SPH algorithms, but do not satisfy all the assumptions of the convergence theorem in \cite{Evers2015}. 

The evolution of an elliptical drop is a benchmark problem for weakly compressible flows, which admits analytical solution \cite{Monaghan1994}. It refers to an initially circular water drop which attains an elliptical shape under a shearing velocity field. The problem involves a conservative part, with $P(\rho) = B \cdot ((\rho/\rho_0)^7 - 1)$, and a dissipative part. The numerical recipe is considered standard in the current SPH literature. For the motion due to the hydrodynamic conservative part we use (\ref{eq: sph_motion1}) and dissipation is modeled with the analogous term of Monaghan and Raffie \cite{Monaghan2013}, which pertains to (\ref{eq: nonlocal_dissipation}). Additionally, the Wendland kernel \cite{Monaghan2013}, and a leapfrog time integrator ---preferred for its symplectic nature--- are used, with $h = 1.5N^{-1/2}$ assigned to all particles. Furthermore, we employ the artificial mass-flux term of Zisis \emph{et al.} \cite{Zisis2015a} with the corresponding parameters $\alpha = 0.5$ and $\beta = 0$, to counteract oscillations in the density profile. We examined two different equations for mass conservation: 1) the discretized temporal evolution of (\ref{eq: mass_conservation}); 2) the discrete SPH-approximation of the continuity equation \cite{Monaghan2013}. Results are practically indistinguishable. The left plot of Figure \ref{fig: droplet} shows the upper half plane of the problem for $N=7232$ at normalized time $t=0.0076$, when Monaghan \cite{Monaghan1994} records the height of the semi-major axis. He finds the latter height $1.91$, compared to $1.95$ of the analytical result (black dashed horizontal line in Figure \ref{fig: droplet}) and our $1.93$. We follow the process until normalized final time $t=0.01$, achieved with a time step $\Delta t= 10^{-6}$. The right plot of Figure \ref{fig: droplet} shows the convergence rates $C^{(2)}_{k+1}$ of the initial measure $\mu_0^N$ and the final one $\mu_T^N$, with respect to the Wasserstein distance between particle systems of successive particle numbers $N_k\in\{2,12,32,52,112,208,448,812,1804,3228,7232\}$. The convergence rates oscillate around the theoretically predicted value $-1/2$, and they tend to become identical. Computing the Wasserstein distance for higher $N_k$ becomes computationally too expensive for our brute-force algorithm. In order to fill the the initial circle (red dashed circle in Figure \ref{fig: droplet}) with $N_k$ particles, we use $\ell_k=\{2,4,6,8,12,16,24,32,48,64,96\}$ particles per unit length to pack particles within a larger square and then disregard all particles falling outside. Recall that the theoretical result (Sections 3 and 4) does not support the current form of dissipation, the equation of state, or the specific functional dependence of $h$ on $N$ ---in the proof $h$ is assumed to be a fixed parameter for all examined $N_k$. Moreover, we obtained indistinguishable results using the SPH-approximation of the continuity equation \cite{Monaghan2013}, for which the theoretical result is not applicable. Yet, our numerical results provide evidence that typically used weakly compressible SPH schemes converge with respect to the Wasserstein distance. Rigorous proofs are left for future work. Our conjecture that the scaling $h\sim N^{-1/d}$ is the correct one, may serve as a guideline.
%
%
%
%
\begin{table}
\centering
\begin{minipage}{0.49\textwidth}
\centering
\includegraphics[width=1.0\textwidth]{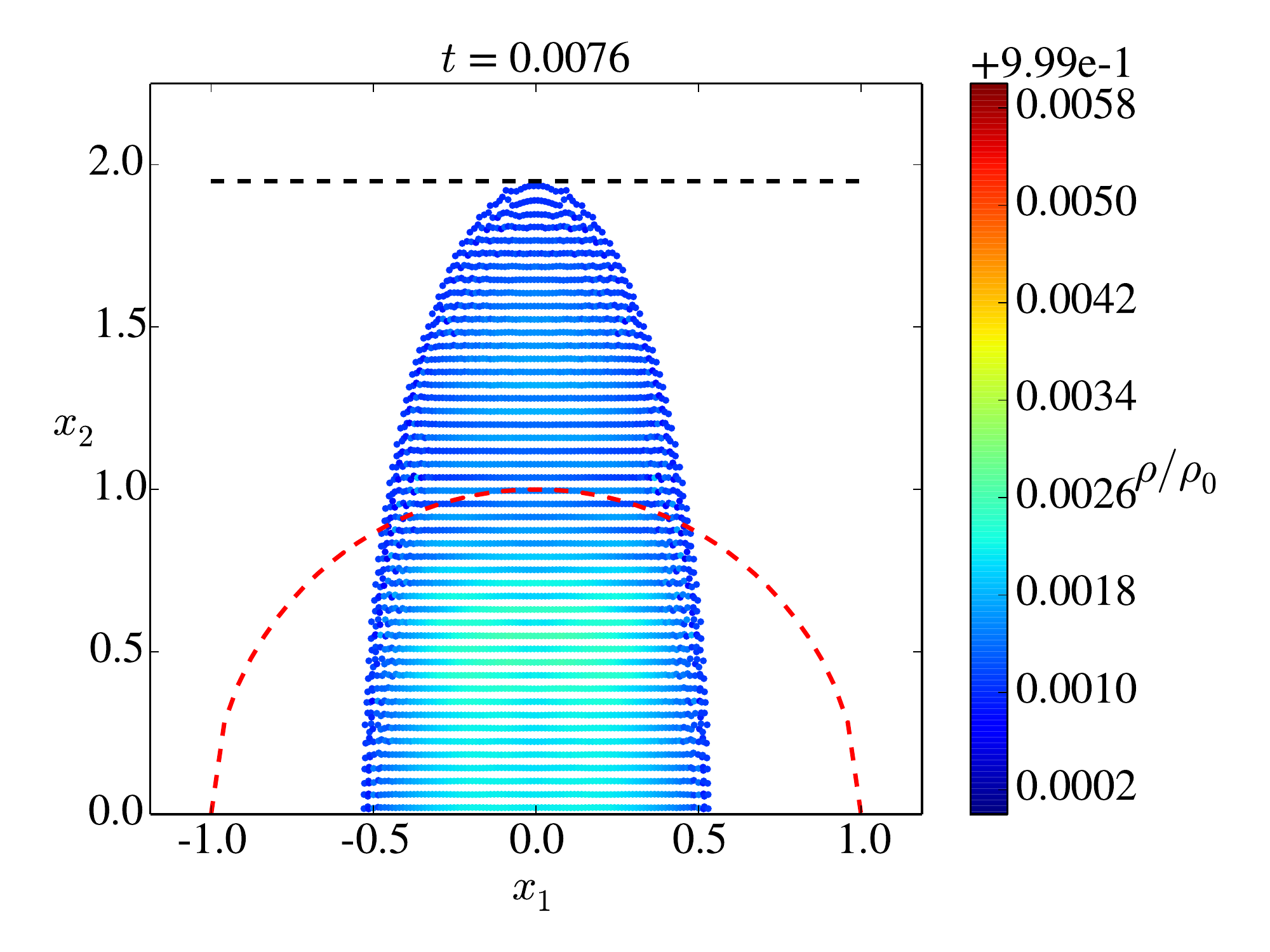}
\end{minipage}
\centering
\begin{minipage}{0.49\textwidth}
\centering
\vspace{-0.15in}
\includegraphics[width=1.0\textwidth]{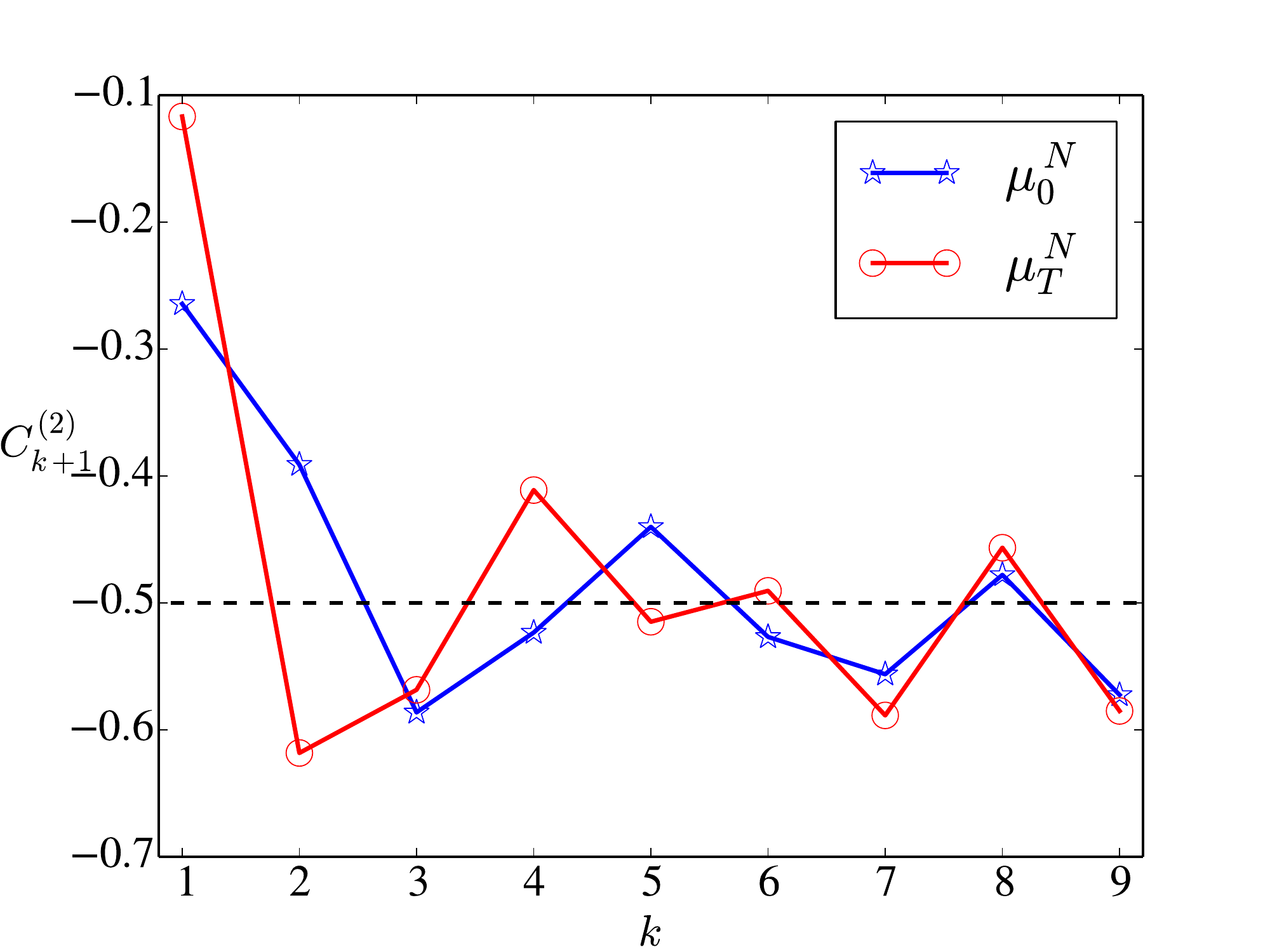}
\end{minipage}
\captionof{figure}{Droplet test; final-to-initial density at t=$0.0076$ and convergence rates of the initial (\emph{blue}) and the final (\emph{red}) measures.}
\label{fig: droplet}
\end{table}
\begin{table}
\centering
\begin{minipage}{0.49\textwidth}
\centering
\includegraphics[width=0.75\textwidth]{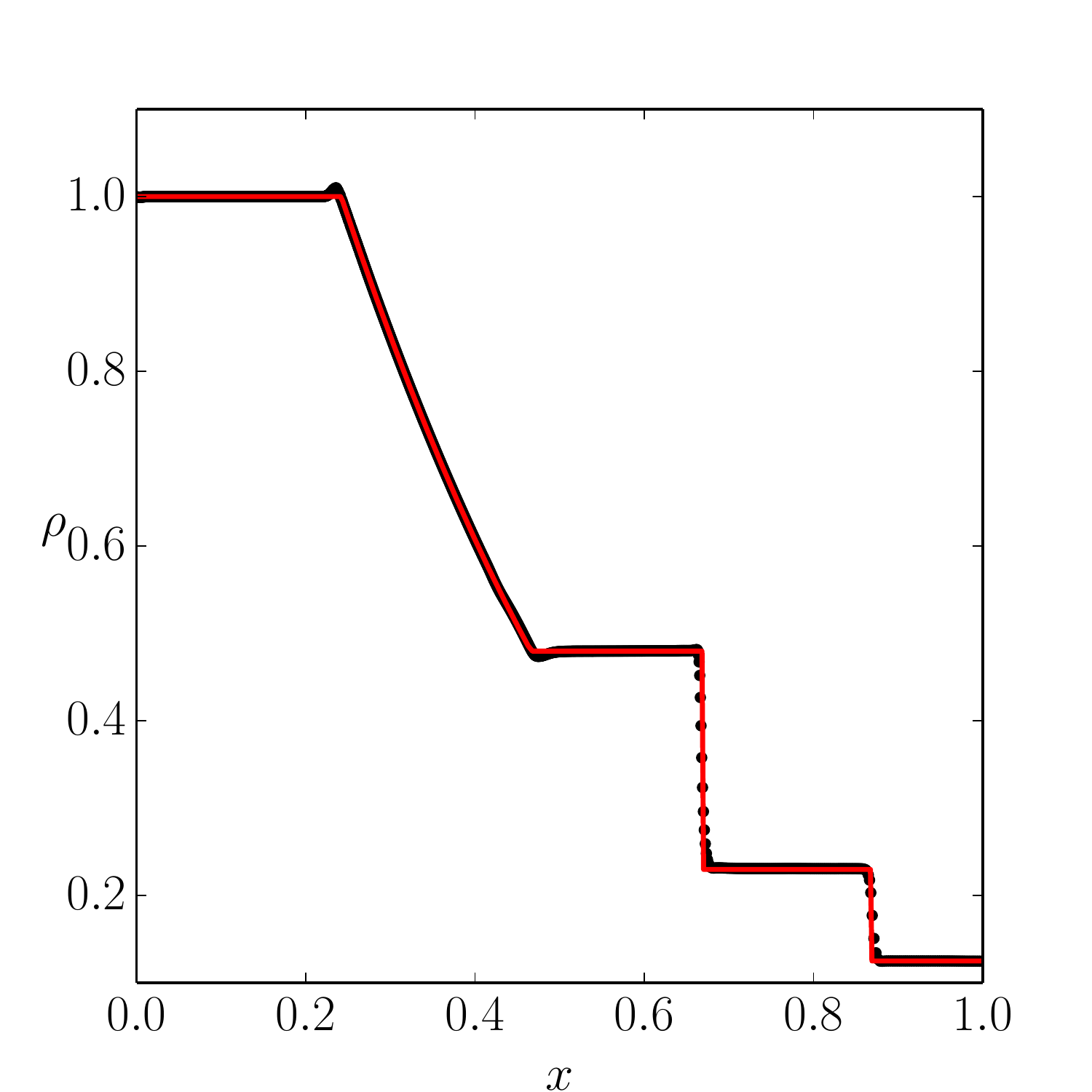}
\end{minipage}
\centering
\begin{minipage}{0.49\textwidth}
\centering
\includegraphics[width=1.0\textwidth]{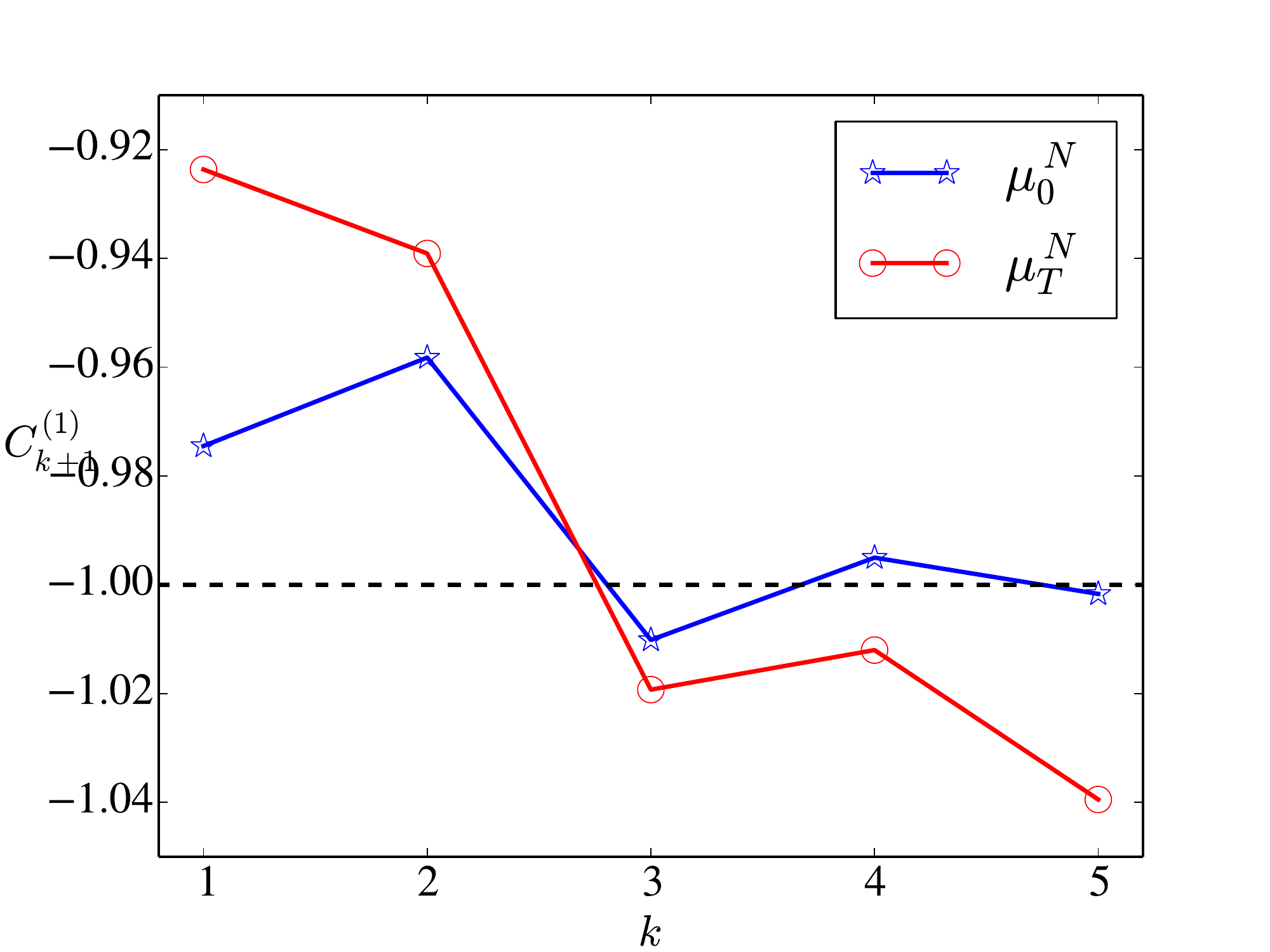}
\end{minipage}
\captionof{figure}{Shock-tube test; density at t=$0.2$ and convergence rates of the initial (\emph{blue}) and the final (\emph{red}) measures.}
\label{fig: shocktube}
\end{table}

The \emph{shock tube test} is a classic one-dimensional test, frequently employed for the validation of fully compressible SPH schemes \cite{Monaghan2005,Price2012,Zisis2015a}. In this test, there is a discontinuity in the density profile of the medium, with $\rho_0(x_0<0.5) = 1$ and $\rho_0(x_0\geq0.5) = 0.125$. We construct the initial density profile with particles of equal masses and solve the SPH system using the differential mass conservation coming from (\ref{eq: mass_conservation}) and the equation of motion (\ref{eq: sph_motion1}), with varying smoothing length $h_i:=1.2 \, m_i/\tilde \rho_i$. The complete solution strategy can be found in Zisis \emph{et al.} \cite{Zisis2015a} and falls within the standard framework  \cite{Monaghan2005,Price2012}. The typical resolution is $450$ particles in total \cite{Price2012,Zisis2015a} and therefore, we examine convergence with respect to the Wasserstein distance for $N_k \in \{18, 45, 90, 225, 450, 900, 1800\}$, to the theoretical value $C^{(1)}_{k+1} = -1$. The distinct characteristic of this case is that $h$ varies spatially and temporally. This is not supported by the theoretical result of convergence \cite{Evers2015}, neither are the ---necessary for the solution--- artificial dissipative terms in all variables. Nonetheless, in Figure \ref{fig: shocktube} the system is shown to converge in a manner similar to the prediction of the theoretical result. The density profile for the highest resolution is also presented in Figure \ref{fig: shocktube}, against the analytical solution (red solid line) at final normalized time $t=0.2$. 

\section{Conclusions}
\label{sec: concl}
The present paper summarizes the authors' results \cite{Evers2015}, regarding the derivation and convergence of SPH. It focuses on describing the three ways to obtain SPH from continuum mechanics via a formulation based on measures. The theoretical convergence is established with respect to the Wasserstein distance, as the number of particles increases, similarly to Di Lisio \emph{et al.} \cite{DiLisio1998}. In fact, the older result is extended by including external fields, local dissipation and non-local interaction forces of the system.

Perhaps the most important limitation of the theoretical proof is that it is not conclusive regarding the convergence of the SPH system as the number of particles grows to infinity and the smoothing length goes to zero at the same time. It rather holds for examining a fixed value of the smoothing length as the number of particles goes to infinity. Additionally, open problems are the inclusion of the following features: dissipation with a non-local character (similar to the one typically used in SPH to mimic viscosity of fluids); equation of states for liquids and solids; spatially varying smoothing length.

Last but not least, the present paper introduces an innovation, by supporting the theoretical proof with numerical evidence. Calculations of the Wasserstein distance reveal that the theoretically predicted convergence rate is observed in the SPH solutions of widely used weakly and fully compressible tests.

\bibliography{SPHbib}

\bibliographystyle{elsarticle-num-names}
\end{document}